\documentclass[a4paper]{spie}
\usepackage[utf8]{inputenc}
\usepackage[round]{natbib}
\usepackage[affil-it]{authblk}
\usepackage{mathrsfs}
\usepackage{amsmath}
\usepackage{amsfonts} %for \mathbb{}
\usepackage{yhmath}
\usepackage{xcolor,colortbl}%for colored table
\definecolor{Gray}{gray}{0.9}%for colored table
\newcolumntype{a}{>{\columncolor{Gray}}c}
\newcolumntype{b}{>{\columncolor{white}}c}
\usepackage{enumitem} %for [label={\alph*)}] and [label={\roman*)}]
\usepackage{booktabs,siunitx} %for \toprule
\usepackage{multirow}
\usepackage{bm}
\usepackage{multicol}
\usepackage{enumerate}
\usepackage{comment}

\title{Experimental performance of graph neural networks on random instances of max-cut}
\author{Weichi Yao,  Afonso S. Bandeira, Soledad Villar}
\begin{document}

\maketitle
\begin{abstract}

%- optimization of hard problems on random inputs
%- GNN appears to be competitive and doesn't require any knowledge of the problem
%- Max-Cut is a great problem to be a test bed because its asymptotics are well understood
%- there are a lot of algorithms that were developed for it that ended up being quite influential (like GW)
%- it is close enough to Community Detection for one to suspect GNN works but different in a very crucial way, that it is an unsupervised problem
%- this way we also investigate the different ways to adapt the framework to do unsupervised things, in a way where we can measure performance cleanly
%- unsupervised ML is the future

This note explores the applicability of unsupervised machine learning techniques towards hard optimization problems on random inputs. In particular we consider Graph Neural Networks (GNNs) -- a class of neural networks designed to learn functions on graphs -- and we apply them to the max-cut problem on random regular graphs. 
We focus on the max-cut problem on random regular graphs because it is a fundamental problem that has been widely studied. In particular, even though there is no known explicit solution to compare the output of our algorithm to, we can leverage the known asymptotics of the optimal max-cut value in order to evaluate the performance of the GNNs. 

In order to put the performance of the GNNs in context, we compare it with the classical semidefinite relaxation approach by Goemans and Williamson~(SDP), and with extremal optimization, which is a local optimization heuristic from the statistical physics literature.  The numerical results we obtain indicate that, surprisingly, Graph Neural Networks attain comparable performance to the Goemans and Williamson SDP. We also observe that extremal optimization consistently outperforms the other two methods. Furthermore, the performances of the three methods present similar patterns, that is, for sparser, and for larger graphs, the size of the found cuts are closer to the asymptotic optimal max-cut value. 

%The Max-Cut %(\textsf{Max-Cut}) problem an NP-hard problem in graph theory that attracted many researchers over the years. Goemans and Williamson proposed their algorithm which combines semi-definite programming and a rounding procedure to produce an approximate solution. Lately, Graph Neural Networks have gained increasing popularity and shown great potential in community detection problems for graph structured data. In the meanwhile, extremal optimization, which is a general-purpose greedy method for approximating solutions to hard optimization problems, has been shown to achieve high accuracy in solving graph bipartitioning problems. 

%In this paper, we adapt the latter two methods for \textsf{Max-Cut} application. In particular, loss functions have been reformulated, using relaxation method and policy gradient method, to allow Graph Neural Networks to work for \textsf{Max-Cut} problem. 

%The cut value produced by the Goemans and Williamson SDP does not differ much from the one produced by Graph Neural Networks. However, there is a a noticeably large difference in the overlap between any of the two methods implies that the three methods always find different cuts, and large fluctuations between individual trials or different initialization for extremal optimization algorithm. 

\end{abstract}
\section{Introduction}

Consider the fundamental problem of \textit{max-cut}, where given a graph, the task is to find a partition of the vertices into two classes such that the number of edges across classes is maximized.
More precisely, given a weighted undirected graph $G=(V,E)$ with $V = \{1, ..., n\}$ and with weights $w_{ij}=w_{ji}\ge 0$, the goal of the max-cut problem is to find a partition $(S, V-S)$ which maximizes the sum of the weights of edges between $S$ and $V-S$.
One can formulate the problem as the following quadratic program:
\begin{align}
    &\max\;\;\frac{1}{2}\sum_{i<j}w_{ij}(1-x_ix_j)\nonumber\\
    &s.t.\;\;\;\;\;x_i\in\{-1,+1\},\;\forall\;i\in[n]\label{Problem}\tag{$P$}
\end{align}
where the $x_i$ are binary variables that indicate set membership, i.e., $x_i = 1$ if $i \in S$ and $x_i = -1$ otherwise.
This integer quadratic program is known to be NP-complete \citep{completeness}.

%Note that the max-cut problem is equivalent to clustering vertices in the graph complement, and so it is not surprising that max-cut enjoys important applications in machine learning.
To date, several algorithms have been developed to approximately solve the max-cut problem, and some approaches, such as the Goemans--Williamson semidefinite relaxation \citep{SDP}, have inspired many algorithms for a variety of problems.
Since optimal solutions of max-cut on sparse graphs tend to correspond to (almost) balanced partitions, this combinatorial optimization problem is closely related to max-bisection, or min-bisection on the complement graph.
In fact, there has been a recent flurry of work on variations of these problems on graph with planted structures.
For example, the stochastic block model is a popular modification of the Erd\H{o}s--R\'{e}nyi random graph, featuring a community structure, in which each edge probability is determined by whether the two vertices reside in a common planted community.
The analysis of max-cut (or min/max -bisection) under such random models of data is convenient because there is a known ground truth to evaluate against, and this feature has contributed to the popularity of this approach \citep{abbe2017community}.
%However, there is another approach available that has received considerably less attention, but has important and fundamental ties to computational complexity, statistical physics, and machine learning.

In this paper, we consider the max-cut problem for random regular graphs without planted structures.
Since there are no planted communities in the underlying random graph distribution, there is no known ground truth to compare against.
Instead, we leverage the well-understood asymptotics of max-cut under this distribution by comparing the max-cut objective value of a given clustering to the asymptotic optimal value~\citep{PA, zdeborova2010conjecture}.
We note that this framework has fundamental ties with various disciplines, and provides a mechanism for identifying computational-to-statistical gaps \citep{bandeira2018notes}.
%In statistical physics, such a setup allows for detailed predictions on the performance of various algorithms.
In machine learning, the absence of a planted signal allows, e.g., one to study false positives in clustering.

Throughout machine learning, it has become popular to solve problems using neural networks, and so we are particularly interested in understanding the performance of such techniques in solving unplanted optimization problems; more specifically, we study the power of \textit{graph} neural networks (GNNs) in solving random max-cut instances.
GNNs are a type of neural network that directly operates on the graph structure. They where first proposed in \citet{GNN1,GNN2}, and have emerged as a powerful class of algorithms to perform complex graph inference leveraging labeled data and tackling classification problems.
Recently, GNNs were used to perform community detection in the stochastic block model, where they proved to be competitive with a variety of other solvers \citep{LGNN}.
%For this reason, one might expect GNNs to similarly perform well in our \textit{unsupervised} setting, in which there are no planted communities.
In this paper, we illustrate that GNNs also perform well in our \textit{unsupervised} setting (despite the lack of labelled data), and we further investigate different ways to adapt the GNN framework so as to perform unsupervised learning.

For the sake of comparison, we solve our max-cut instances using GNNs, as well as the Goemans--Williamson SDP approach and extremal optimization.
Extremal optimization (EO) is a local method for combinatorial optimization problems that first appeared in the statistical physics literature~\citet{EO1}.
It has been successfully applied to solve the graph bi-partitioning problem \citep{EOGBP} and the traveling salesman problem \citep{chen2007optimization}, among other discrete optimization problems.
In the case of the bi-partitioning problem on random graphs, simulations performed in \cite{EOGBP} suggest that EO converges to near optimal configurations in $O(N)$ time where $N$ is the number of nodes, but approximation guarantees for this algorithm are not currently known.

Another algorithmic approach worth mentioning is a message passing algorithm developed by Montanari for optimization of the Sherrington-Kirkpatrick Hamiltonian~\citep{montopt}. The optimization of the Sherrington-Kirkpatrick (SK) Hamiltonian can be viewed as a ``gaussian'' analogue of max-cut on random graphs; and it was its study that allowed for the asymptotic understanding of max-cut on random $d$-regular graphs~\citep{PA}. This algorithm is proved~\citep{montopt} (conditioned on a mild statistical physics conjecture) to be asymptotic optimal in the SK setting. While it would be fascinating to understand the performance of an analogue of the algorithm in~\cite{montopt} to the max-cut problem in random sparse graphs, this is outside of the scope of this paper.

The paper is organized as follows.
In Section 2, we describe the Goemans and Williamson algorithm, the extremal optimization algorithm, and how we adapt the Line Graph Neural Network model from \citet{LGNN} to the max-cut problem.
Section 3 shows extensive numerical results, where the performance of the three methods is evaluated and compared.
In Section 4, we conclude with an outlook on future work.

\section{Algorithms for Max-Cut problems}
%Let $G = (V, E)$ be a random graph with vertex set $|V| = n$, and let $W \in \{0, 1\}^{n\times n}$ denote its adjacency matrix. 
\subsection{The Goemans-Williamson Algorithm}
In their very influential paper, \citet{SDP} introduce an approximation algorithm for the max-cut problem based on semidefinite programming. The algorithm consists of the following steps:
\\(i) Relax the optimization problem \eqref{Problem} to the vector programming \eqref{ProblemRlx}: 
    \begin{align}
        &\max_{\mathbf u}\;\;\frac{1}{2}\sum_{i<j}w_{ij}(1-\mathbf{u}_i^T\mathbf{u}_j)\nonumber\\
        &s.t.\;\;\;\;\;\left\Vert\mathbf{u}_i\right\Vert^2=1,\;\text{and }\mathbf{u}_i\in \mathbb{R}^n,\;\forall\;i\in[n]\tag{$P^\prime$}\label{ProblemRlx}
    \end{align}
    which is equivalent to the following semidefinite program \eqref{ProblemSDP}:
    \begin{align}
        &\max_X\;\;\frac{1}{2}\sum_{i<j}w_{ij}(1-X_{ij})\nonumber\\
        &s.t.\;\;\;\;\;X\succeq 0,\;\text{and }X_{ii}=1,\;\forall\;i\in[n].\tag{SDP}\label{ProblemSDP}
    \end{align}
    Solve \eqref{ProblemRlx} and denote its solution as $\mathbf{u}_i$, $i\in[n]$. \\
     (ii) Let $r\in\mathbb{R}^n$ be a vector uniformly distributed on the unit sphere \label{rounding1}
    \\
    (iii) Produce a cut subset $S = \{i|r^T\mathbf{u}_i\ge 0\}$. \label{rounding2}

%Note that one can run the randomized rounding procedure (ii)-(iii) several times and select the best cut.
\citet{SDP} showed that their algorithm is able to obtain a cut whose expected value is guaranteed to be no smaller than a particular constant $\alpha_{GW}$ times the optimum cut:
\begin{align*}
    \textsf{max-cut}(G)\ge\alpha_{\mathrm{GW}}\dfrac{1}{2}\sum_{i<j}w_{ij}(1-\mathbf{u}_i^T\mathbf{u}_j)\ge\alpha_{\mathrm{GW}}\textsf{max-cut}(G)
\end{align*}
with $\alpha_{\mathrm{GW}} = \min_{-1\le x\le 1}\frac{\frac{1}{\pi}\arccos(x)}{\frac{1}{2}(1-x)}\simeq 0.878$. %The constant $\alpha_{GW}$ is referred to as the approximation ratio.

\subsection{Extremal Optimization}
\citet{EO1,EO2} introduced a local search procedure called extremal optimization (EO), modeled after the Bak-Sneppen mechanism \citep{BSmodel}, which was originally proposed to describe the dynamics of co-evolving species.% According to the authors, the goal of EO is to provide a fast and reliable approximation method for finding optimal or  near-optimal solutions with high probability when optimizing a system of many variables with respect to some cost function, which often exhibits a complex ``landscape'' %See Landscape Paradigms in Physics and Biology, edited by H. Frauenfelder, et al. (Elsevier,Amsterdam, 1997) from https://arxiv.org/pdf/cond-mat/0104214.pdf
%in configuration space, in cases where the relation between individual components of the system is frustrated \citep{CostFunct}.
%It is a method based on the dynamics of non-equilibrium processes and in particular those exhibiting self-organized criticality \citep{selforganized}, where better solutions emerge dynamically without the need for parameter tuning.

In \citet{EOGBP} the EO algorithm is applied to the graph bi-partitioning problem (GBP) and it is empirically suggested that EO efficiently approximates the optimal bi-partition on Erd\H{o}s--R\'{e}nyi and random 3-regular graphs. Specifically, what their experiments actually show, is that the EO produces a partition which cut value scales in the same way than the expected value of the optimal cut with the number of nodes. They also remark that the convergence of EO appears to be linear on the number of steps. However, mathematical guarantees for the EO algorithm are not yet known.  %EO, like other local search methods [26], gets stuck ever more often in poor local minima. However, when we consider the best out of multiple runs with the EO algorithm, we recover results close to those predicted by theoretical arguments.

It is straightforward to adapt the EO algorithm from GBP to max-cut. In statistical physics, max-cut and GBP can be formulated in terms of finding the ground state configurations of an Ising model on random $d$-regular graphs. For any graph, the Ising spin Hamiltonian is given by 
\begin{align*}
    \mathcal{H}(\mathbf{x}) = -\sum_{(i,j)\in E}J_{ij}x_ix_j
\end{align*}
where each spin variable $x_i\in\{-1,+1\}$ is connected to each of its nearest neighbors $j$ via a bond variable $J_{ij}\in\{-1,+1\}$, assigned at random. The configuration space $\Omega$ consist of all configurations $\mathbf{x}=(x_1,x_2,...,x_n)\in\Omega$ where $\left\vert\Omega \right\vert=2^n$. The goal of EO is to minimize the cost function $ \mathcal{H}(\mathbf{x})$.

Note that,
\begin{align*}
    \textsf{max-cut}(G)=\max_{\mathbf{x}\in\{\pm 1\}^n}\frac{1}{2}\sum_{i<j}w_{ij}(1-x_ix_j)=\left\vert E(G)\right\vert-\min_{\mathbf{x}\in\{\pm 1\}^n}\frac{1}{2}\sum_{i<j}w_{ij}x_ix_j.
\end{align*}
Instead of finding a ground state of the ferromagnetic Ising model with magnetization fixed at zero as in the GBP, the max-cut problem is to find a ground state configurations of the anti-ferromagnetic Ising model by minimizing the Hamiltonian $\mathcal{H}$ that can be written as 
\begin{align*}
    \mathcal{H}(\mathbf{x}) = \frac{1}{2}\sum_{i<j}w_{ij}x_ix_j
\end{align*}
where $w_{i,j}$ is the $(i,j)$ entry of the adjacency matrix $W=W^T$ of the graph. To find low-energy configurations, EO assigns a fitness to each $x_i$
\begin{align*}
    \lambda_i = -x_i\left(\frac{1}{2}\sum_jw_{ij}x_j\right),
\end{align*}
so that
\begin{align*}
   \mathcal{H}(\mathbf{x}) = -\sum_{i=1}^n\lambda_i.
\end{align*}
heuristically, the fitness of each variable assesses its contribution to the total cost. Note that $\lambda_i=-\frac{d}{2}+b_i$, where $g_i$ is the number of edges connecting $i$ to other vertices within its same set (``good'' edges), and $b_i$ is the number of and edges connecting $i$ to vertices across the partition (``bad'' edges). We consider a normalized version, 
\begin{align}\label{eq:normlambda}
    \lambda_i = \frac{b_i}{g_i+b_i}=\frac{b_i}{d}
\end{align}
as in \citet{EOGBP}, so that $\lambda_i \in [0,1]$. 

EO proceeds through a local search of $\Omega$ by sequentially changing variables with ``bad'' fitness on each update. It ranks all the variables $x_i$ according to fitness $\lambda_i$, i.e., find a permutation $\Pi$ of the variable labels $i$ with
\begin{align*}
    \lambda_{\Pi(1)}\le\lambda_{\Pi(2)}\le\cdots\le\lambda_{\Pi(n)}.
\end{align*}
If one always only updates the lowest-ranked variable, the search will risk reaching a ``dead end'' in the form of a poor local minimum. In contrast, the idea is to consider a probability distribution over the ranks $k$,
\begin{align}\label{eq:tau}
    P_k\propto k^{-\tau},\;\;\;1\le k\le n, 
\end{align}
for a given value of the parameter $\tau$. At each update, select a rank $k$ according to $P_k$, then if advantageous for fitness, change the state of the variable $x_i$ with $i=\Pi(k)$. After each update, the fitness of the changed variable and of all its neighbours are re-evaluated according to (\ref{eq:normlambda}). %Repeat the fitness evaluation and update as long as desired.

\subsection{Graph Neural Network}
Recently, \citet{LGNN} proposed a family of Graph Neural Networks (GNNs) \citep{GNN1,GNN2} for community detection on graphs. Inspired by the success of regularized spectral methods based on non-backtracking random walks for community detection by~\citet{krzakala2013spectral}, ~\cite{saade2014spectral} and \citet{bordenave2015non}, \citet{LGNN} extend the Graph Neural Networks architecture to express such objects.  
%GNNs to operate on the line graph using the non-backtracking operator (LGNN) and have shown that such GNNs can reach the statistical and computational signal-to-noise detection thresholds for random graph families such as the stochastic block model (SBM), 
These GNNs learn clustering algorithms in a purely (supervised) data-driven manner without access to the underlying generative models, and even empirically improving upon the state of the art for some specific regimes. They have also been adapted to problems like graph matching and traveling salesman~\cite{nowak2018revised}.

The generic GNN, introduced in \citet{GNN0} and later simplified in \citet{GNN3,GNN4,GNN5} is a neural network architecture that is based on finding a way to combine local operators on a undirected graph $G = (V, E)$ to provide the desired output. The usual choice of local operators are degree operator $(Dx)_i:= \mathrm{deg}(i) \cdot x_i$, $D(x) = \mathrm{diag}(A\mathbf{1})x$, and the adjacency operator $A$, which is the linear map given by the adjacency matrix $A_{i,j} = 1$ when $(i, j) \in E$. \citet{LGNN} further consider the power graph adjacency  $A_j = \min(1, A^{2^j})$, which encodes $2^j$-hop neighborhoods into a binary graph and therefore allows to combine and aggregate local
information at different scales. The architecture has two folds, one is the operation on the original graph $G$, the other is on the line graph $L(G) = (V_L, E_L)$, which represents the edge adjacency structure of $G$, with non-backtracking operator $B$. The vertices $V_L$ of $L(G)$ are the ordered edges in $E$, that is $V_L = \{(i \rightarrow j); (i, j) \in E\} \bigcup \{(j \rightarrow i); (i, j) \in E\}$, so $|V_L| = 2|E|$. The non-backtracking operator $B \in \mathbb{R}^{2|E|\times 2|E|}$ encodes the edge adjacency structure as follows,
\begin{equation*}
B_{(i\rightarrow j),(i^\prime\rightarrow j^\prime)}=
   \begin{cases}
   1 &\text{if }j = i^\prime\text{ and } j^\prime\neq i,\\
   0 &\text{otherwise.}
   \end{cases}
\end{equation*}
with the corresponding degree $D_B = \mathrm{diag}(B\mathbf{1})$ operators. This effectively defines edge features that are diffused and updated according to the edge adjacency of $G$. Edge and node features are combined at each layer using the edge indicator matrices Pm, Pd $\in\{0, 1\}^{|V|\times 2|E|}$, defined as Pm$_{i,(i\rightarrow j)} = 1$, Pm$_{j,(i\rightarrow j)} = 1$, Pd$_{i,(i\rightarrow j)} = 1$, Pd$_{j,(i\rightarrow j)} = -1$ and 0 otherwise. Then, for each layer, with an input signal $u^{(k)} \in \mathbb{R}^{|V|\times b_k}$ on $G$ and $v^{(k)} \in \mathbb{R}^{|V_L|\times b_k}$ on line graph $L(G)$, the model produces $u^{(k+1)} \in\mathbb{R}^{|V|\times b_{k+1}}$ and $v^{(k+1)} \in\mathbb{R}^{|V_L|\times b_{k+1}}$ as
\begin{align}\label{eq:LGNNmodel}
    u^{(k+1)}_{i,l} &= \rho\left[u_i^{(k)}\theta_{1,l}^{(k)}+\sum_{j=0}^{J-1}(A^{2^j}u^{(k)})_i\theta_{3+j,l}^{(k)}+ \{\mathrm{Pm,Pd}\}v^{(k)}\theta^{(k)}_{3+J,l}\right],l=1,...,b_{k+1}/2,i\in V,\nonumber\\
    v^{(k+1)}_{i^\prime,l}&=\rho\left[v_{i^\prime}\gamma_{1,l}^{(k)}+(D_{L(G)}v^{(k)})_{i^\prime}\gamma_{2,l}^{(k)}+\sum_{j=0}^J(A_{L(G)}^{2^{j}}v^{(k)})_{i^\prime}\gamma_{3+j,l}^{(k)}+[\{\mathrm{Pm,Pd}\}^Tu^{(k+1)}]_{i^\prime}\gamma_{3+J,l}^{(k)}\right],i\in V_L,\nonumber\\
    u^{(k+1)}_{i,l} &=u_i^{(k)}\theta_{1,l}^{(k)}+\sum_{j=0}^{J-1}(A^{2^j}u^{(k)})_i\theta_{3+j,l}^{(k)}+\{\mathrm{Pm,Pd}\}v^{(k)}\theta^{(k)}_{3+J,l},l=b_{k+1}/2+1,...,b_{k+1},i\in V,\nonumber\\
    v^{(k+1)}_{i^\prime,l}&=v_{i^\prime}\gamma_{1,l}^{(k)}+(D_{L(G)}v^{(k)})_{i^\prime}\gamma_{2,l}^{(k)}+\sum_{j=0}^J(A_{L(G)}^{2^{j}}v^{(k)})_{i^\prime}\gamma_{3+j,l}^{(k)}+[\{\mathrm{Pm,Pd}\}^Tu^{(k+1)}]_{i^\prime}\gamma_{3+J,l}^{(k)},i\in V_L 
\end{align}
where $\Theta= \{\theta_1^(k), . . . , \theta_{J+3}^{(k)},\gamma^{(k)}_1,...,\gamma^{(k)}_{J+3}\}$, $\theta^{(k)}_s$, $\gamma^{(k)}_s \in \mathbb{R}^{b_k\times b_{k+1}}$ are trainable parameters and $\rho(\cdot)$ is a point-wise nonlinearity, chosen in this work to be $\rho(z) = \max(0, z)$. The model \eqref{eq:LGNNmodel} may seem a little daunting, but it's basically an unrolled overparametrized power iteration, on node features as well as edge features, where the parameters $\Theta$ indicate how to combine the operators of the graph to produce an optimized \emph{regularized spectral method}.  

Assume that a training set $\{(G_t, \mathbf{y}_t)\}_{t\le T}$ is given, with any signal $\mathbf{y} = (y_{1},...,y_{n}):V\rightarrow \{1,2,...,C\}^n$ encoding a partition of $V$ into $C$ groups. It is then used by LGNN to learn a model $\widehat{\mathbf{y}}=\Phi(G,\Theta)$ trained by minimizing
\begin{align} \label{loss}
    L(\Theta) = \frac{1}{T}\sum_{t\le T}l(\Phi(G_t,\Theta),\mathbf{y}_t).
\end{align}
In the max-cut problem, $C=2$. Unlike community detection on SBM, the max-cut problem on random graphs is cast as an unsupervised learning problem with no knowledge of $\{\mathbf{y}_t\}$ in the training process. We define the loss function as the cut itself, then the model $\widehat{\mathbf{y}}=\Phi(G,\theta)$ is learned by training 
\begin{align*}
    \max_{\Theta} L(\Theta) = \max_{\Theta}\frac{1}{T}\sum_{t\le T}l(\Phi(G_t,\Theta))=-\min_\Theta \frac{1}{T}\sum_{t\le T}-l(\Phi(G_t,\Theta)).
\end{align*}
with 
\begin{align}\label{eq:Lundif}
    l(\Phi(G_t,\Theta)) = \frac{1}{4}\mathbf{x}^T\mathcal{L}_G\mathbf{x}=:f(\mathbf{x},G_t)
\end{align}
where $\mathbf{x} = (x_1,x_2,...,x_n)\in\{\pm 1\}^n$ is the resulted configurations and where $\mathcal{L}_G$ is the graph Laplacian. Note that, $f(\mathbf{x},G_t)$ is not differentiable, but we wish to apply stochastic gradient decent methods. To this end, we propose two different methods: a relaxation method, and a policy gradient method. We make use of the probability matrix $\pi(\mathbf{x}|\Theta,G)\in [0,1]^{n\times 2}$ (with $(i,j)$ entry equal to $\mathbb{P}(x_i=j|\Theta,G)$), $j=\pm 1$, which is an intermediate result produced by LGNN.

\subsubsection{Relaxation method on loss function for \textsf{Max-Cut} problem}
%In \textsf{Max-Cut} problem, $C=2$. 
Denote $p_i = %\mathbb{P}(x_i=-1|\Theta,G)$ or $
\mathbb{P}(x_i=1|\Theta,G)$. Substitute $\mathbf{x}$ with $2\mathbf{p}-1\in [-1,1]^n$, then we can write the loss function as
\begin{align*}
    l(\Phi(G,\Theta)) =\frac{1}{4}(2\mathbf{p}-1)^T\mathcal{L}_G(2\mathbf{p}-1)
\end{align*}
which is differentiable with respect to the parameters of the model.
\subsubsection{Policy gradient methods on loss function for \textsf{Max-Cut} problem}
%\asb{I find this section difficult to read}
Policy gradient is a technique from reinforcement learning that relies upon optimizing parameterized policies with respect to the expected return by gradient descent. Informally the ``policies" we consider in this context correspond with the distribution of choices of cuts given a graph, which is a function of the set of parameters $\Theta$. The expected return under that policy is the expected value of the cut given the graph. 

In the max-cut problem, for any given graph $G$, consider the following optimization
\begin{align*}
    \Theta^\ast = \arg\max_\Theta\mathbb{E}_{\mathbf{x}\sim\pi_{\Theta}(\cdot|G)}[f(\mathbf{x},G)].
\end{align*}
%with 
%\begin{align*}
%    \mathbb{E}_{\mathbf{x}\sim\pi_{\Theta}(\cdot|G)}[f(\mathbf{x},G)]=\int f(\mathbf{x},G)\pi_{\Theta}(\mathbf{x}|G)d\mathbf{x}.
%\end{align*}
Here $\mathbf{x}\sim\pi_{\Theta}(\cdot|G)$ can be viewed as the policy and $\mathbb{E}_{\mathbf{x}\sim\pi_{\Theta}(\cdot|G)}[f(\mathbf{x},G)]$ as the expected return.

Many machine learning approaches base their optimization on variations of policy gradients (see for instance \cite{PG1,pg2}). We use the simplest policy gradient formulation, corresponding to the seminal work by \cite{PG}. Note that 
\begin{align}
\nabla_\Theta\mathbb{E}_{\mathbf{x}\sim\pi_{\Theta}(\cdot|G)}[f(\mathbf{x},G)]&=\int f(\mathbf{x}, G)\pi_{\Theta}(\mathbf{x}|G)\nabla_\Theta\log\pi_{\Theta}(\mathbf{x}|G)d\mathbf{x}\nonumber\\
&=\mathbb{E}_{\mathbf{x}\sim\pi_{\Theta}(\cdot|G)}[f(\mathbf{x},G)\nabla_\Theta\log\pi_{\Theta}(\mathbf{x}|G)]\label{eq:policy}.
\end{align}

The intermediate output probability matrix $\pi(\cdot|\Theta,G)$ from LGNN enables us to sample $\mathbf{x}$ from its distribution (that depends on the current $\Theta$ and the given training graph $G$). We sample $K$ independent instances of $\mathbf x$  obtaining $\{\mathbf{x}_k\}_{k=1}^K$. Then (\ref{eq:policy}) can be approximated by
\begin{align}
\nabla_\Theta\mathbb{E}_{\mathbf{x}\sim\pi_{\Theta}(\cdot|G)}[f(\mathbf{x},G)]
&\approx\dfrac{1}{K}\sum_{k=1}^K[f(\mathbf{x}_k,G)\nabla_\Theta\log\pi_\Theta(\mathbf{x}_k|G)]\\
&=\nabla_\Theta\left\{\dfrac{1}{K}\sum_{k=1}^K[f(\mathbf{x}_k,G)\log\pi_\Theta(\mathbf{x}_k|G)]\right\}. \label{gradient}
\end{align}
The loss function can then set to be
\begin{align*}
l(\Phi(G,\Theta)) = \frac{1}{K}\sum_{k=1}^K[f(\mathbf{x}_k,G)\log\pi_\Theta(\mathbf{x}_k|G)],
\end{align*}
with $\{\mathbf{x}_k\}_{k=1}^K$ sampled from $\pi_\Theta(\cdot|G)$; and we can use the gradient from \eqref{gradient} to implement stochastic gradient descent on the loss function \eqref{loss}. 

\section{Experimental evaluation}
%Random regular graphs provide the simplest model on which this intuition can be made precise~\asb{What is ``this intution'' refering to here?}. 
Denoting by $G
^{\text{reg}}(n, d)$ the uniform distribution over graphs with $n$ vertices and uniform degree $d$, \citet{PA} have proved that, with high probability, as $n\rightarrow \infty$ the size of the max-cut satisfies 
\begin{align*}
    \textsf{MaxCut}(G^{\text{Reg}}(n,d)) =n\left(\frac{d}{4} + P_\ast\sqrt{\frac{d}{4}}+o_d\left(\sqrt{d}\right)\right)+o(n)
\end{align*}
% \begin{align*}
%     \dfrac{\textsf{MaxCut}(G^{\text{Reg}}(n,d))}{n} = \dfrac{d}{4}+P_\ast\sqrt{\frac{d}{4}} + o_d\left(\sqrt{d}\right)
% \end{align*}
with $P_\ast \approx 0.7632$ an universal constant related to the ground state energy of the Sherrington-Kirkpatrick model, that can be expressed analytically via Parisi’s formula~\citep{PA,percus2008peculiar, talagrand2006parisi}.

Motivated by these asymptotics, for any candidate solution $\mathbf{x} = (x_1,...,x_n)$, and the corresponding cut size $z=\frac{1}{2}\sum_{i<j}w_{ij}(1-x_ix_j)$, we evaluate it by computing its corresponding $P$ as 
\begin{align}\label{eq:P}
    P = \dfrac{z/n - d/4}{\sqrt{d/4}}.
\end{align}
The larger the $P$ the better the cut, and moreover, if $P$ is closer $P_\ast \approx 0.7632$ we expect the candidate cut to be among the best possible in the graph.

It is worth noting that the SDP relaxation is not tight in the asymptotic limit \citep{montsdp}. Actually, in the limit the corresponding $P$ of the SDP's (fractional) solution is 1, which is not an improvement over a simple spectral method. However our experiments in the next Section indicate that, at least for the parameter range we are investigating, the rounding procedure still produces non-trivial cuts.

We are also interested in how much the resulted configurations for a given graph overlap with each other. We define the overlap measure $\nu$ for any two configuration $\mathbf{x}_i=(x_1,x_2,...,x_n)\in \{\pm 1\}^n$, $i=1,2$, as 
\begin{align}\label{eq:overlap}
    \nu = \frac{1}{n}\left\vert \langle\mathbf{x}_1,\mathbf{x}_2\rangle \right\vert.
\end{align}
%For the resulted configurations from any of the two methods, the mean value of the overlap as well as the standard deviation on the test graphs will be reported and evaluated.

%\subsection{Simulation results}
In our numerical simulations, we mainly focus on the performance of the discussed methods when the graphs are relatively sparse. To obtain a more thorough comparison of the performance of the three methods, various values of node size $n$ and degree $d$ were chosen; the cut values were computed over 1000 random $d$-regular graphs with $n$ vertices. 

The values of the tuning parameters of the three methods used in the simulations are given and explained as follows. We use SDPNAL+~\citep{sdpnal} to solve the semidefinite program~ \eqref{ProblemSDP}, and the random rounding procedure chooses the best cut over 500 iterations. For EO, there is only one parameter, the exponent $\tau$ in the probability distribution in equation (\ref{eq:tau}), governing the update process and consequently the performance of EO. \citet{EOGBP} investigated numerically the optimal value for $\tau$, which they view as a function of runtime $t_{\max}$ and node size $n$. They observed optimal performance on graphs of node size $n\le 10000$ with $t_{\max}\ge 100n$ with $\tau \approx 1.3 \sim 1.6$. Based on their numerical findings, they claim that the value of $\tau$ does not significantly change the value of the solutions as long as the algorithm runs for enough iterations. Therefore, following their discussion, we simply fixed $\tau = 1.4$ and set a sufficiently large number of iterations $t_{\max} = 10^4n$. In addition, as the labels are randomly assigned at the beginning of the EO algorithm, it is expected that the different initialization may lead to different outcomes. Therefore, multiple runs of the EO is suggested. Here we chose the best run out of two. To train \citet{LGNN} model on the max-cut problem, we set the input signals to be $x^{(0)} = \mathrm{deg}(G)$ and $y^{(0)} = \mathrm{deg}(L(G))$. Other parameters were set to be the same as what was used in \citet{LGNN}'s simulations, $J = 3$, $K=30$, $b_k=10$, $k = 2,...,K-1$ ($b_0=1$ and $b_K = C$, where $C=2$ for the max-cut problem) in (\ref{eq:LGNNmodel}). We trained the model on 5000 graphs with the same node size $n$ and degree $d$. Though different sets of tuning parameters in GNNs may change the overall performance, we believe such settings still can provide us a reasonable understanding of how well GNNs perform compared to the other two methods. 

All our code is publicly available in \cite{code}.

\subsection{The size of the computed cut}
\paragraph{Same node size $n$ with different degrees $d$.}
First we explore the setting where we fix the node size $n$, and let the degree $d$ change. Table \ref{tab:sameNdiffd} provides the results on graphs with $n = 500$ nodes. The results on graphs of other node sizes are similar.

\begin{table}[h!]
    \centering
    \begin{tabular}{cccccc}
         \toprule
         \multirow{3}{*}{$n$}&\multirow{3}{*}{$d$} & \multicolumn{4}{c}{Methods}  \\ 
         \cmidrule{3-6}
         &&\multicolumn{2}{c}{GNNs} & \multirow{2}{*}{SDP} & \multirow{2}{*}{EO} \\
         \cmidrule{3-4}
         && Relaxation & Policy Gradient & & \\
         \midrule
         \multirow{5}{*}{500}
         &20 &0.6136	&0.6259	&0.6742	&0.7315\\
         \cmidrule{2-6}
         &15 &0.6070	&0.6289	&0.6783	&0.7359\\
         \cmidrule{2-6}
         &10 &0.6704	&0.5989	&0.6820	&0.7352\\
         \cmidrule{2-6}
         &5 &0.7014	    &0.6682	&0.6898	&0.7369\\
         \cmidrule{2-6}
         &3 & 0.7074	&0.6928	&0.7015	&0.7266\\
         \bottomrule
    \end{tabular}
    \caption{Computed $P$ values for different methods, with same node size $n = 500$ and various degrees $d$.}
    \label{tab:sameNdiffd}
\end{table}

The experiments reported in Table \ref{tab:sameNdiffd} show that the computed $P$ improves as the degree decreases. %, i.e., all methods appear to result in a better cut (when normalized  when the graphs get sparser. 
 Overall, EO appears to give the best performance. SDP outperforms GNNs while GNNs seem to catch up when the graph gets sparser.

\paragraph{Same degree $d$ with increasing node size $n$}
Asymptotically, the optimal $P$ on random $d$-regular graphs should approach $P_\ast$ as $n\rightarrow\infty$ and $d \rightarrow \infty$. Tables \ref{tab:samed3diffN} and \ref{tab:samed10diffN} provide two instances when $d=3$ and $d = 10$, respectively. In both cases, the computed $P$ for all three methods increase as the node size $n$ increases, getting closer to $P_\ast$. Furthermore, the EO method shows greater advantage over the other methods as the node size $n$ increases. The performance of SDP is still comparable to that of GNNs when node size $n$ is up to 500.
\begin{table}[h!]
    \centering
    \begin{tabular}{cccccc}
         \toprule
         \multirow{3}{*}{$d$}&\multirow{3}{*}{$n$} & \multicolumn{4}{c}{Methods}  \\ 
         \cmidrule{3-6}
         &&\multicolumn{2}{c}{GNN} & \multirow{2}{*}{SDP} & \multirow{2}{*}{EO} \\
         \cmidrule{3-4}
         && Relaxation & Policy Gradient & & \\
         \midrule
         \multirow{4}{*}{3}
         &50 & 0.6716&	0.6868&	0.6981& 0.6985\\
         \cmidrule{2-6}
         &100 &0.6978&	0.6914&	0.7090&	0.7118\\
         \cmidrule{2-6}
         &200 &0.7075&	0.7010&	0.7091&	0.7210\\
         \cmidrule{2-6}
         &500 &0.7074&	0.6928&	0.7015&	0.7266\\
         \bottomrule
    \end{tabular}
    \caption{Computed $P$ values for different methods, with same degrees $d=3$ but different node sizes $n$.}
    \label{tab:samed3diffN}
\end{table}

\begin{table}[h!]
    \centering
    \begin{tabular}{cccccc}
         \toprule
         \multirow{3}{*}{$d$}&\multirow{3}{*}{$n$} & \multicolumn{4}{c}{Methods}  \\ 
         \cmidrule{3-6}
         &&\multicolumn{2}{c}{GNN} & \multirow{2}{*}{SDP} & \multirow{2}{*}{EO} \\
         \cmidrule{3-4}
         && Relaxation & Policy Gradient & & \\
         \midrule
         \multirow{4}{*}{10}
         &50  &0.5874&	0.5995&	0.6614&	0.6643\\
         \cmidrule{2-6}
         &100 &0.6296&	0.6368&	0.6889&	0.7033\\
         \cmidrule{2-6}
         &200 &0.6456&	0.6598&	0.6919&	0.7241\\
         \cmidrule{2-6}
         &500 &0.6704&	0.5989&	0.6820&	0.7369\\
         \bottomrule
    \end{tabular}
    \caption{Computed $P$ values for different methods, with same degrees $d=10$ but different node sizes $n$.}
    \label{tab:samed10diffN}
\end{table}

\section{Conclusions}

The experimental results presented in this paper suggest that unsupervised machine learning techniques can be successfully adapted to hard optimization problems on random inputs. We observe this in the particular case of max-cut on random regular graphs, where graph neural networks learn algorithms that attain comparable performance to the Goemans and Williamson SDP for this particular model. This finding is noteworthy due to the unsupervised --and computationally hard-- nature of the task.

Overall, the three methods we consider --SDP, EO, and GNNs--  appear to compute solutions with comparable objective values. Preliminary numerical simulations, not reported in the current paper but available in~\cite{code}, suggest that the overlap (or correlation) between different solutions with similar objective value can still be quite different, which is a statement on the difficulty of the optimization problem. %We also observe that the performances of the three methods present a similar pattern, that is, for sparser, and for larger graphs, the size of the found cuts are closer to the asymptotic optimal max-cut value.

Throughout these simulations we observe that extremal optimization is the most computationally efficient and it consistently outperforms the SDP and GNNs under both, dense and sparse regimes. We therefore believe there is a need for further theoretical analysis of extremal optimization. In our simulations EO appeared to be somewhat initialization dependent, which if not a finite $n$ artifact, could prove a difficulty when trying to theoretically analyze this method.

Finally, there exists significant room for improvement on the implementation of the GNNs. It would not be surprising if a better architecture or choice of hyperparameters can be used to improve upon our results. 

\section*{Acknowledgments}
The authors would like to thank Dustin Mixon and Zhengdao Chen. SV is partly supported by NSF-DMS 1913134, EOARD FA9550-18-1-7007 and the Simons Algorithms and Geometry (A\&G) Think Tank. ASB was partially supported by NSF grants DMS-1712730 and DMS-1719545, and by a grant from the Sloan
Foundation.

\bibliographystyle{unsrtnat}
\bibliography{ref}
\end{document}